\documentclass[10pt,leqno]{amsart}
\usepackage{graphicx}
\usepackage{indentfirst,csquotes}

\usepackage{amssymb,amsthm,amsmath}
\usepackage{xcolor,paralist,hyperref,titlesec,fancyhdr,etoolbox}
\newtheorem{theorem}{Theorem}[]

\newtheorem{proposition}{Proposition}

\titleformat{\section}[display]{\normalfont\huge\bfseries\centering}{\centering\chaptertitlename\thechapter}{10pt}{\Large}
\titlespacing*{\section}{0pt}{0ex}{0ex}

\hypersetup{ colorlinks=true, linkcolor=black, filecolor=black, urlcolor=black }

\usepackage{lipsum}

\begin{document}
\title[Inverse source  problems  for a multidimensional  time - fractional wave . . .]{Inverse source  problems  for a multidimensional  time - fractional wave   equation with integral overdetermination conditions} 
\author[D.K.Durdiev]{D.K.Durdiev}
\date{\today}
\address{$ ^1$Bukhara Branch of Romanovskii Institute of Mathematics,
Uzbekistan Academy of Sciences, Bukhara, Uzbekistan,\\
$ ^2$Bukhara State University, Bukhara, 705018 Uzbekistan}
\email{d.durdiev@mathinst.uz, durdiev65@mail.ru}
\maketitle

\let\thefootnote\relax

\begin{abstract}
In this paper, we consider two linear inverse problems for the time-fractional wave equation, assuming that its right-hand side takes the separable form \( f(t)h(x) \), where \( t \geq 0 \) and \( x \in \Omega \subset \mathbb{R}^N \). The objective is to determine the unknown function \( f(t) \) (Inverse Problem 1) and \( h(x) \) (Inverse Problem 2), given that the other function is known.

For Inverse Problem 1, we impose an overdetermination condition in the form of a spatial integral over the domain \( \Omega \), involving the solution of the corresponding direct problem—an initial-boundary value problem with standard Cauchy conditions and homogeneous Dirichlet boundary conditions. The integral is weighted by the known spatial factor \( h(x) \) from the right-hand side of the equation. This choice of an additional condition enables the explicit construction of a solution to the inverse problem and allows us to prove its unique solvability within the class of regular solutions. To study the direct problem, the regular solution approach is used.

For Inverse Problem 2, we introduce a novel integral-type additional condition, referred to as the time-averaged velocity, incorporating an appropriate weight function. The time-dependent factor of the right-hand side of the equation serves as this weight function. Depending on its choice, the additional condition reduces to specifying either the final-time offset or the time-averaged offset. Under this formulation, we establish a new uniqueness result.
\end{abstract} 

\bigskip


\begin{center}
    {\bf Introduction}
\end{center}

Many real-world problems in natural sciences cannot be accurately modeled using differential equations with standard integer-order derivatives. However, replacing these derivatives with fractional-order counterparts often leads to more precise and realistic descriptions. In this context, fractional operators, particularly the Riemann-Liouville and Caputo derivatives, play a crucial role in mathematical modeling.

Fractional diffusion and diffusion-wave equations have been widely used to describe anomalous diffusion phenomena, such as superdiffusion and subdiffusion, which are not well captured by classical diffusion models. These equations have gained significant attention in recent years \cite{Dur01}-\cite{Dur03}. In particular, time-fractional diffusion and wave equations with variable coefficients have attracted considerable interest due to their applications in various scientific and engineering disciplines, including physics, chemistry, biology, control theory, and viscoelasticity (see \cite{Dur1}-\cite{Dur3} and references therein).

The direct problems, i.e., initial value problem and initial-boundary value problems
for the time fractional diffusion and wave equations have been studied extensively in recent years, for instance, on maximum principle
\cite{Dur04}, on some uniqueness,  existence results and  analytic solutions \cite{Kil} (see also the references therein).

In this paper, we examine the initial-boundary value problem for the time-fractional diffusion-wave equation:
\begin{equation}\label{eq(10.1)}
  \partial_t^{\alpha}u-Lu=f(t)h(x), \, \, \,
  \alpha\in(1, 2),\, \, \, \, (t, x)\in \Omega_T:=(0, T]\times\Omega, \end{equation}
\begin{equation}\label{eq(10.2)}
  u(0, x)=\varphi(x), \, \, \, \, \,  u_t(0, x)=\psi(x), \, \, \, \, \,  x\in\overline{\Omega},
 \end{equation}
\begin{equation}\label{eq(10.3)}
   u(t, x)=0, \, \, \,  (t, x)\in \partial\Omega_T:=[0, T]\times\partial\Omega,
\end{equation}
where $L=div\left[\mathbf{A}(x) \nabla u(t, x)\right]-c(x)u(t, x)$ with $\mathbf{A}(x)=\left(a_{ij}(x)\right)_{i, j=1, \cdots, N},$ $\mathbf{A}=\mathbf{A}^{\top}$ (the symbol $\top$ denotes the transpose),
   \, \, $\sum\limits_{i, j=1}^N a_{ij}(x)\xi_i\xi_j\geq\delta\sum\limits_{i=1}^N\xi_i^2,$ \,\, $\delta=const>0,$ \, $x\in \Omega,$ \, $\xi\in \mathbb{R}^N,$
 in addition, it is
assumed $a_{ij}(x), \, c(x)$ satisfy some conditions of smoothness and $c(x)\geq 0;$
   \, $\Omega\subset\mathbb{R}^N$ is a bounded
domain with smooth boundary $\partial\Omega,$ \, $N\geq 1$;  the Caputo fractional differential operator
$\partial_t^{\alpha}$ of the order $1<\alpha<2$  is defined in \cite[pp. 90-99]{Kil}:
\begin{equation*}\partial_t^{\alpha}u(x, t):=I_{0+}^{2-\alpha}u_{tt}(x, t)=\frac{1}{\Gamma(2-\alpha)}\int_0^t\frac{u_{\tau\tau}(x, \tau)}{(t-\tau)^{\alpha-1}}d\tau,\end{equation*}
\begin{equation*}I_{0+}^{\mu}u(x, t):=\frac{1}{\Gamma(\mu)}\int_0^t\frac{u(x, \tau)}{(t-\tau)^{1-\mu}}d\tau, \, \, \mu\in(0, 1),\end{equation*}
the function $\Gamma(\cdot)$ represents the Gamma function, while $I_{0+}^{\mu}u(x, t)$ denotes the Riemann–Liouville fractional integral of the function $u(x, t)$ with respect to the variable $t.$

To study the problem (\ref{eq(10.1)})-(\ref{eq(10.3)}), we employ  {\it the regular solution approach} and introduce classes of continuous and smooth functions to represent the solution. Let $C(\Omega_T),$ $C\left(\overline{\Omega}_T\right)$ be spaces of continuous functions on $\Omega_T,$ $\left(\overline{\Omega}_T\right),$ respectively.
Introduce also  the classes of functions \\
$$C^{1,\, 0 }_{t,\, x}\left(\overline{\Omega}_T\right):=\Big\{u(t, x)\in C\left(\overline{\Omega}_T\right): \,  u_t(t, x)\in  C\left(\overline{\Omega}_T\right) \Big\},$$
$$C^{\alpha,\, 2}_{t, \, x}\left(\Omega_T\right):=\Big\{u(t, x)\in C^{1,\, 0 }_{t, \, x}\left(\overline{\Omega}_T\right): \, \partial_t^{\alpha}u, \, \, u_{xx}\in  C\left(\Omega_T\right) \Big\}.$$

Given  $\varphi(x),$ $\psi(x),$ $f(t), $ $h(x),$ $\alpha$ and $L$ we  call (\ref{eq(10.1)})-(\ref{eq(10.3)}) as a {\bf Direct Problem (DP)}  and we will look for its solution in  $C^{\alpha,\, 2}_{t, \, x}\left(\Omega_T\right)\cap C^{1,\, 0 }_{t,\, x}\left(\overline{\Omega}_T\right).$

For this problem, we  investigate the following inverse problems:

\textbf{Inverse Problem 1 (IP1):} Given the operator $L$, the initial states $\varphi,$ $\psi$, the number $\alpha\in(1, 2)$ and  the source
function $h$,  determine a pair of functions $\left\{u(t, x), \, f(t)\right\},$ satisfying the conditions $u(t, x)\in C^{\alpha,\, 2}_{t, \, x}\left(\Omega_T\right)\cap C^{1,\, 0 }_{t,\, x}\left(\overline{\Omega}_T\right),$ \, $f(t)\in C[0, T],$ equations (\ref{eq(10.1)})-(\ref{eq(10.3)}), as well as the overdetermination condition
\begin{equation}\label{eq(ad1)}
 \int\limits_{\Omega} h(x)u(t, x)dx=g(t), \, \, \, \, 0\leq t \leq T,
\end{equation}
where  $g(t)$ is given  function.
This type of measurement represents the
 average of $u$ with weight function $h$ over the spatial  domain $\Omega$. And there are many literatures on inverse problems
by this type measured data for fractional diffusion and diffusion-wave equations (see \cite{Kil11}-\cite{Kil13}
and the references therein).

\textbf{Inverse Problem 2 (IP2):}  Given the operator $L$, the initial states $\varphi,$ $\psi$, the number $\alpha\in(1, 2)$ and  the source
function $f$,  find a pair of functions $\left\{u(t, x), \, h(x)\right\},$ satisfying the conditions $u(t, x)\in C^{\alpha,\, 2}_{t, \, x}\left(\Omega_T\right)\cap C^{1,\, 0 }_{t,\, x}\left(\overline{\Omega}_T\right),$ \, $h(x)\in C(\overline{\Omega}),$
equations (\ref{eq(10.1)})-(\ref{eq(10.3)}) and the equality
\begin{equation}\label{eq(adc3)}
\int\limits_0^T f(t)u_t(t, x)dt=\omega(x), \, \, \, \, x\in \overline{\Omega},
\end{equation}
where  $\omega(x)$ is a known  function.
In this problem, the integral conditions (\ref{eq(adc3)}) is an  additional condition for determining the function $h(x).$ The non-standard overdetermination  condition (\ref{eq(adc3)}) corresponds to the time-averaged velocity with the weight function $f(t).$ If $f(t)\in H^1(0, T)$ ($H^1(0, T)$ is the Sobolev space of functions which consists of all integrable real functions $f: (0, T)\mapsto \mathbb{R}$ such that $f, \, f'\in L^2(0, T)$) and
$f(T)\neq 0,$ then integrating by parts,  (\ref{eq(adc3)}) can be rewritten as
\begin{equation}\label{eq(adcc3)}
f(T)u(T, x)-\int\limits_0^T f'(t)u(t, x)dt=\omega(x)+f(0)\varphi(x), \, \, \, \, x\in \overline{\Omega}.
\end{equation}
Thus, specifying an overdetermination condition of the form (\ref{eq(adc3)}) is equivalent to specifying a special combination of the conditions of the final   time-measured output  and the time-averaged output (the integral on the left side of (\ref{eq(adcc3)})). Note that condition (\ref{eq(adc3)}) for $f(t)=const$ coincides with specifying the final offset in time with the  output $\omega(x)+const\cdot\varphi(x)$ and for $f(t)\in H^1_0(0, T)$ the offset averaged in time with the weight $- f'(t)$.

An inverse source problem for equation (\ref{eq(10.1)}) involves determining the unknown functions \( f(t) \) and/or \( h(x) \) based on additional measurement data related to \( u \). These problems can generally be classified into inverse \( x \)-source problems and inverse \( t \)-source problems, depending on whether the spatial or temporal component of the source needs to be reconstructed.

Numerous studies have addressed inverse \( x \)-source and \( t \)-source problems for time-fractional diffusion-wave equations of the form (\ref{eq(10.1)}) with \( N=1 \) and \( L=\frac{\partial^2}{\partial x^2} \), considering various types of fractional derivatives and overdetermination conditions. Most existing works on inverse problems for fractional equations focus on the Riemann-Liouville and Caputo derivatives with orders \( \alpha \in (0,1) \). A comprehensive review of theoretical results and computational approaches for inverse source problems in time-fractional diffusion equations with \( \alpha \in (0,1) \) can be found in \cite{Rundell, Rundell1}.
In \cite{Yamo}, the stability of an inverse source problem for a fractional diffusion-wave equation with \( \alpha \in (0,2] \) and a symmetric, uniformly elliptic operator with \( x \)-independent smooth coefficients is analyzed. Using eigenfunction expansions, the authors establish stability results for recovering the time-dependent factor of the source from single-point observations over the time interval \( (0,T) \). In \cite{Yamo1}, inverse source problems for a broad class of parabolic, hyperbolic, and time-fractional evolution equations are studied under partial interior observations. By assuming that the source term is separable in space and time, the authors prove uniqueness for simultaneous reconstruction of both components without requiring nonzero initial conditions. Furthermore,
the works \cite{Yamo2, Yamo3} focus on the inverse problem of identifying the spatially dependent source term in a time-fractional diffusion-wave equation using additional subboundary measurement data. Uniqueness is demonstrated through analyticity arguments, asymptotic behavior of solutions, and a newly established unique continuation principle for time-independent coefficients. The uniqueness of determining a spatially varying factor from decay conditions as \( t \to \infty \) is proven in \cite{Yamo2}, assuming the source is inactive during the observation period. In \cite{Yamo3}, Lipschitz stability is derived using an appropriate function space topology based on the adjoint system of the fractional diffusion-wave equation. Moreover, in \cite{Yamo4}, an inverse problem for a subdiffusion equation with a Riemann-Liouville fractional derivative and a general elliptic operator in an arbitrary multidimensional domain is studied. Using Fourier decomposition, the authors prove existence and uniqueness theorems for the classical solution of the initial-boundary value problem and the unique determination of an \( x \)-dependent source function. The work \cite{Kil122} investigates the uniqueness of an inverse source problem for linear time-fractional diffusion equations with time-dependent coefficients, where the unknown spatially varying source function is reconstructed from final-time measurement data.

Several studies have also explored inverse problems for time-dependent sources. In \cite{Dur7}, the inverse problem of determining a time-dependent source term in a one-dimensional time-fractional diffusion equation from energy measurements is examined. The work \cite{Dur8} extends this analysis to multi-dimensional settings using boundary Cauchy data. In \cite{Dur9}, the recovery of a time-dependent coefficient in a generalized time-fractional diffusion equation with non-local boundary and integral overdetermination conditions is considered. The existence and uniqueness of the solution are established using function expansions in a bi-orthogonal system and fractional calculus techniques. The work \cite{Dur10} addresses forward and inverse source problems for a two-dimensional time-fractional diffusion equation, where the inverse problem is formulated as an Abel integral equation of the first kind, later transformed into a second-kind integral equation using fractional differentiation.

In this paper, we begin by investigating an inverse \( t \)-source problem for equation (\ref{eq(10.1)}) using integral measurement data over the spatial domain. We first analyze the unique solvability of the direct problem (DP) using Fourier decomposition method. To solve inverse problem IP1, we employ the theory of integral equations. Unlike previous works, the integral overdetermination condition (\ref{eq(ad1)}) (where \( h(x) \) represents the spatial component of the applied load \( F(t,x)=f(t)h(x) \)) allows us to construct an explicit solution to IP1 and prove its well-posedness in the class of regular solutions. Furthermore, using techniques based on a priori estimates for functions with fractional derivatives, we establish uniqueness results for inverse problem IP2, which involves determining the spatially dependent source function (\( x \)-source problem). To the best of our knowledge, no prior studies have considered inverse problems IP1 and IP2 under the specific additional conditions formulated in this work for determining \( f(t) \) and \( h(x) \), respectively.

The structure of the paper is as follows. Section 2 provides some preliminaries concerning with fractional calculus.  In Section 3, we analyze the direct problem using the Fourier decomposition method. Section 4 focuses on the existence and uniqueness of solutions to IP1, employing fractional calculus and integral equation techniques. Although uniqueness is typically established before existence in inverse problem theory, we first derive key results from integral equation theory in Section 4 and use them in Section 5 to prove the uniqueness of the solution to IP1. Section 6 investigates uniqueness results for IP2, where a space-dependent source is determined. Finally, Section 7 provides concluding remarks.

\begin{center}
\textbf{Preliminaries.}
\end{center}
The Mittag-Leffler function $E_{\alpha,
\beta}(z)$ is  defined by the following series \cite[pp. 40-45]{Kil}:
$$
E_{\alpha,
\beta}(z):=\sum_{n=0}^{\infty}\frac{z^n}{\Gamma(\alpha{n}+\beta)},$$
  $\Gamma(\cdot)$ is the  Euler's gamma function $\alpha, \, z, \,  \beta\in\mathbb{C}$, $
\mathfrak{R}(\alpha)>0.$
\begin{proposition}\label{prop1}
\cite[pp. 40-45]{Kil} Let $0<\alpha<2$ and $\beta\in\mathbb{R}$ be arbitrary. We suppose that $\kappa$ is such that $\pi\alpha/2<\kappa<\min\{\pi,\pi\alpha\}$. Then there exists a constant $C=C(\alpha,\beta,\kappa)>0$ such that
$$
\left|E_{\alpha,\beta}(z)\right|\leq\frac{C}{1+|z|},\quad \kappa\leq|\mbox{arg}(z)|\leq\pi.
$$
\end{proposition}
\begin{proposition}\label{prop2}
The following differentiation formulas are valid:
$$\frac{d}{dt}E_{\alpha, 1} \left(-\lambda_n t^{\alpha}\right)=-\lambda_n t^{\alpha-1} E_{\alpha, \alpha} \left(-\lambda_n t^{\alpha}\right),$$
$$\frac{d}{dt}\left[tE_{\alpha, 2} \left(-\lambda_n t^{\alpha}\right)\right]= E_{\alpha, 1} \left(-\lambda_n t^{\alpha}\right),$$
$$\frac{d}{dt}\left[t^{\alpha-1}E_{\alpha, \alpha} \left(-\lambda_n t^{\alpha}\right)\right]= t^{\alpha-2}E_{\alpha, \alpha-1} \left(-\lambda_n t^{\alpha}\right).$$
\end{proposition}

{\it Proof.}  This proof follows from the formula 1.10.7 in \cite[p. 50]{Kil}.
\begin{proposition}\label{prop3}
The following fractional differentiation formulas are valid:
$$\partial_t^{\alpha-1}\left[t^{\alpha-1}E_{\alpha, \, \alpha}\left(-\lambda_n t^{\alpha}\right)\right]=-E_{\alpha, \, 1}\left(-\lambda_n t^{\alpha}\right),$$
  $$\partial_t^{\alpha-1}E_{\alpha, \, 1}\left(-\lambda_n t^{\alpha}\right)=-\lambda_ntE_{\alpha, \, 2}\left(-\lambda_n t^{\alpha}\right).$$
\end{proposition}

{\it Proof.} Based on 2.1.54 in  \cite[p. 78]{Kil}, for $\alpha-1\in (0, 1)$  we have
$$\frac{\partial^{\alpha-1}}{\partial t^{\alpha-1}}\left[t^{\alpha-1}E_{\alpha, \, \alpha}\left(-\lambda_n t^{\alpha}\right)\right]=E_{\alpha, \, 1}\left(-\lambda_n t^{\alpha}\right),$$
where $\frac{\partial^{\alpha-1}}{\partial t^{\alpha-1}}y(t)$ is the Riemann--Liouville fractional
derivative  of a  suitable function $y(t):$
$$ \frac{\partial^{\alpha-1}}{\partial t^{\alpha-1}}y(t)=\frac{1}{\Gamma(2-\alpha)}\frac{d}{d  t}\int_0^t(t-\tau)^{1-\alpha}y(\tau)d\tau.$$
The Caputo fractional derivative for  $0<\alpha-1<1$  with the Riemann-Liouville fractional derivative
 is connected by
\begin{equation}\label{FD}
\partial_t^{\alpha-1}y(t)=\frac{\partial^{\alpha-1}}{\partial t^{\alpha-1}}y(t)-\frac{y(+0)}{\Gamma(2-\alpha)t^{\alpha-1}}, \, \, \, t>0.
\end{equation}
Since $y(t)=t^{\alpha-1}E_{\alpha, \, \alpha}\left(-\lambda_n t^{\alpha}\right)$ and $y(+0)=0,$ then, in this case $\partial_t^{\alpha-1}y(t)=\frac{\partial^{\alpha-1}}{\partial t^{\alpha-1}}y(t).$ This yields the first formula of Proposition \ref{prop3}.

By the same argument (see the  formula 2.1.54 in  \cite[p. 78]{Kil}),  we get
$$\frac{\partial^{\alpha-1}}{\partial t^{\alpha-1}}E_{\alpha, \, 1}\left(-\lambda_n t^{\alpha}\right)=t^{1-\alpha}E_{\alpha, \, 2-\alpha}\left(-\lambda_n t^{\alpha}\right).$$
Next, we convert the right side of this equality as follows
$$t^{1-\alpha}E_{\alpha, \, 2-\alpha}\left(-\lambda_n t^{\alpha}\right)=t^{1-\alpha}\sum_{n=0}^{\infty}\frac{\left(-\lambda_n t^{\alpha}\right)^n}{\Gamma(\alpha{n}+2-\alpha)}=\frac{t^{1-\alpha}}{\Gamma(2-\alpha)}-\lambda_n tE_{\alpha, \, 2}\left(-\lambda_n t^{\alpha}\right).$$
If we take into account $E_{\alpha, \, 1}(+0)=1$  and (\ref{FD}), then from the last relations follows the validity of the second equality of Proposition \ref{prop3}.

\begin{proposition}\label{prop4}
 For $1<\alpha<2$ and $\lambda>0,$ if $f(t)\in AC[0, T],$ then, we have
$$ \partial_t^{\alpha-1}\int_0^tf(s)(t-s)^{\alpha-2}E_{\alpha, \, \alpha-1}\left(-\lambda_n (t-s)^{\alpha}\right)ds\quad\quad\quad\quad\quad\quad\quad\quad\quad\quad\quad\quad
\quad\quad\quad\quad\quad$$$$=f(t)-\lambda\int\limits_0^t f(s)(t-s)^{\alpha-1}E_{\alpha, \, \alpha}\left(-\lambda (t-s)^{\alpha}\right)d s, \quad 0<t\leq T.$$
\end{proposition}

{\it Proof.} We carry out the proof similarly to Lemma 2.9 in the paper \cite{Kil1}. Let
$$F(t):=\int_0^tf(s)(t-s)^{\alpha-2}E_{\alpha, \, \alpha-1}\left(-\lambda_n (t-s)^{\alpha}\right)ds.$$
By Proposition \ref{prop1}, we know $\left|E_{\alpha, \ \alpha-1}(-\eta)\right|\leq C$ for $\eta>0.$ Further, we have $|F(t)|\leq(C/(\alpha-1))t^{\alpha-1}\|f\|_{C[0, T]}.$ Then, we define $F(0)=0,$   it is easy to show that $F(t)\in C[0, T]$. Integrating by part, we have
$$F(t)=-f(0)t^{\alpha-1}E_{\alpha, \, \alpha}\left(-\lambda_n t^{\alpha}\right)-\int_0^tf'(s)(t-s)^{\alpha-1}E_{\alpha, \, \alpha}\left(-\lambda_n (t-s)^{\alpha}\right)ds.$$
Further, by $f(t)\in AC[0, T]$ and Proposition \ref{prop2},  we calculate
$$F'(t)=-f(0)t^{\alpha-2}E_{\alpha, \, \alpha-1}\left(-\lambda_n t^{\alpha}\right)-\int_0^tf'(s)(t-s)^{\alpha-2}E_{\alpha, \, \alpha-1}\left(-\lambda_n (t-s)^{\alpha}\right)ds\in L(0, T),$$
therefore  $F(t)\in AC[0, T].$ We note that for function $F(t)$ with $F(0)=0$ the Caputo derivative $ \partial_t^{\alpha-1}$ coincides with Riemann-Liouville derivative $\frac{\partial^{\alpha-1}}{\partial t^{\alpha-1}}$ of the order ${\alpha-1}\in(0, 1)$ (see (\ref{FD})).
That is
\begin{equation}\label{eq(44.4)}\partial_t^{\alpha-1}F(t)=\frac{\partial^{\alpha-1}}{\partial t^{\alpha-1}}F(t)=\frac{1}{\Gamma(2-\alpha)}\frac{d}{dt}\int\limits_0^t\frac{F(\tau)}{(t-\tau)^{\alpha-1}}.
\end{equation}

Next, we find
$$\int\limits_0^t\frac{F(\tau)}{(t-\tau)^{\alpha-1}}=\int\limits_0^t\frac{d\tau}{(t-\tau)^{\alpha-1}}\int\limits_0^{\tau}f(s)E_{\alpha, \, \alpha-1}\left(-\lambda_n (\tau-s)^{\alpha}\right)\frac{ds}{(\tau-s)^{2-\alpha}}$$
$$\quad\quad\quad\quad\quad=\int\limits_0^tf(s)ds\int\limits_s^{t}\frac{E_{\alpha, \, \alpha-1}\left(-\lambda_n (\tau-s)^{\alpha}\right)d\tau}{(t-\tau)^{\alpha-1}(\tau-s)^{2-\alpha}}$$
$$\quad\quad\quad\quad\quad=\int\limits_0^tf(s)ds\sum\limits_{k=o}^{\infty}\frac{(-\lambda)^k}{\Gamma(\alpha k+\alpha-1)}\int\limits_s^{t}\frac{d\tau}{(t-\tau)^{\alpha-1}(\tau-s)^{2-\alpha-\alpha k}}.$$
It is not difficult to note
$$\int\limits_s^{t}\frac{d\tau}{(t-\tau)^{\alpha-1}(\tau-s)^{2-\alpha-\alpha k}}=(t-s)^{\alpha k}\int\limits_0^{1}\tau^{-2+\alpha+\alpha k}(1-\tau)^{1-\alpha}d\tau$$
$$=(t-s)^{\alpha k}B(\alpha k+\alpha-1, \ 2-\alpha)=(t-s)^{\alpha k}\frac{\Gamma(\alpha k+\alpha-1)\Gamma(2-\alpha)}{\Gamma(\alpha k+1)},$$
where $B(\cdot, \cdot)$ is the  Euler's beta function.

With this in mind, we  continue
$$\int\limits_0^t\frac{F(\tau)}{(t-\tau)^{\alpha-1}}=\Gamma(2-\alpha)\int\limits_0^tf(s)\sum\limits_{k=o}^{\infty}\frac{\left(-\lambda(t-s)^{\alpha}\right)^k}{{\Gamma(\alpha k+1)}}ds$$\quad\quad\quad\quad\quad\quad\quad\quad\quad\quad\quad\quad\quad\quad\quad\quad\quad\quad\quad\quad\quad\quad\quad$$=\Gamma(2-\alpha)\int\limits_0^tf(s)E_{\alpha, \, 1}\left(-\lambda_n (t-s)^{\alpha}\right)ds.$$
Thus, from (\ref{eq(44.4)}) based on Proposition \ref{prop2},
we get
$$\partial_t^{\alpha-1}F(t)=f(t)-\lambda\int\limits_0^tf(s)(t-s)^{\alpha-1}E_{\alpha, \, \alpha}\left(-\lambda (t-s)^{\alpha}\right)ds.$$
Proposition \ref{prop4} is proven.

$AC^n[0, T]$ represents the class of functions $v(x)$ that are continuously differentiable on $[0, T]$ up to order $n-1$, such that $v^{(n-1)}(x) \in AC[0, T]$. Here, $n = 1, 2, \dots$, and $AC[0, T]$ denotes the class of absolutely continuous functions on $[0, T]$

\begin{proposition}\label{prop5}
 For any function $\vartheta(t)\in AC^2[0, T],$ one has the
inequality
\begin{equation}\label{eq(41.444)}
\int\limits_0^T\vartheta'(t)\partial_t^{\alpha}\vartheta(t)dt\geq\frac{1}{2\Gamma(1-\gamma)}\int\limits_0^T\frac{\vartheta'^2(t)}{(T-t)^{\gamma}}dt-\frac{T^{1-\gamma}}{2\Gamma(\gamma)}\vartheta'^2(0),
 \end{equation} {\it where} $\alpha\in(1, 2),$  \ $\gamma=\alpha-1\in(0, 1).$
\end{proposition}

  {\bf Proof.} First we use the inequality (see \cite{Ali})
  \begin{equation}\label{eq(1.4444)}\vartheta'(t)\partial_t^{\alpha}\vartheta(t)=\vartheta'(t)\partial_t^{\gamma}\vartheta'(t)\geq \frac{1}{2}\partial_t^{\gamma}\vartheta'^2(t).\end{equation}
  Integrating this relation over interval $(0, T)$ in $t$ and transform the integral appearing on the right side as follows
  \begin{equation*}\frac{1}{2}\int\limits_0^T\partial_t^{\gamma}\vartheta'^2(t)dt=\frac{1}{2\Gamma(1-\gamma)}\int\limits_0^T\int\limits_0^t\frac{1}{(t-\tau)^{\gamma}}\frac{\partial \vartheta'^2(\tau)}{\partial \tau} d\tau dt
\end{equation*}
\begin{equation*}\label{eq(1.41)}
=\frac{1}{2\Gamma(1-\gamma)}\int\limits_0^T\frac{\partial \vartheta'^2(\tau)}{\partial \tau}\int\limits_{\tau}^T\frac{dt}{(t-\tau)^{\gamma}}d\tau
=\frac{1}{2(1-\gamma)\Gamma(1-\gamma)}\int\limits_0^T\frac{\partial \vartheta'^2(\tau)}{\partial \tau}(T-\tau)^{1-\gamma}d\tau.
 \end{equation*}
Integrating by parts the last integral yields
$$\int\limits_0^T\partial_t^{\gamma}\vartheta'^2(t)dt=\frac{1}{2\Gamma(1-\gamma)}\int\limits_0^T\frac{\vartheta'^2(t)}{(T-t)^{\gamma}}dt-\frac{T^{1-\gamma}}{2\Gamma(\gamma)}\vartheta'^2(0).$$
Proposition \ref{prop5} is proven.

\begin{center}
\textbf{Study of DP}
\end{center}

Assume that the first derivatives of the functions $a_{ij}(x)$ are H\"{o}lder continuous, and that the function $c(x)$ is also H\"{o}lder continuous in $\Omega$. Now, consider the following spectral problem:\begin{equation}\label{5}
 LX=-\lambda X,\,\,\,x\in \Omega, \quad X\big|_{\partial\Omega}=0.
\end{equation}

It is known (see, \cite{Il30}) that problem $(\ref{5}),$ as non-trivial solutions has orthonormal eigenfunctions $X_{n}(x), \ n \geq1$, which are complete in $ L_2(\Omega).$ The corresponding eigenvalues $\lambda_n$ are positive, arranged in non-decreasing order, and repeated according to their finite multiplicity:
$0<\lambda_1\leq\lambda_2\leq\cdots$.

 Using the  Fourier decomposition method, we construct   the solution to DP (\ref{eq(10.1)})-(\ref{eq(10.3)}). For this aim,  we write
\begin{equation}\label{eq(1.2)}
u(x, t)=\sum_{n=1}^{\infty} u_n(t) X_n(x),
\end{equation}
where $u_n(t)=\left(u(t, x), \ X_n(x)\right):=\int_{\Omega} u(t, x) X_n(x) d x$ is the scalar product of functions $u(t, \cdot)$ and $X_n(\cdot)$ in $ L_2(\Omega).$

Evidently, for each $n=1, 2, \cdots,$ $u_n(t)$ is the solution to the following  ordinary fractional differential equation:
\begin{equation}\label{eq(ef3)}
\partial_t^{\alpha}u_n(t)+\lambda_n u_n(t)=f(t)h_n,
\end{equation}
where   $h_n= \left( h,  X_n\right). $ \ In accordance with  (\ref{eq(10.2)})  the solution of equation (\ref{eq(ef3)}) must satisfy  the initial conditions
\begin{equation*}
u_n(0)=\varphi_n=\left( \varphi,  X_n\right), \quad\quad u'_n(0)=\psi_n=\left(\psi,  X_n\right).
\end{equation*}

The solution of this problem, for each $n=1, 2, \cdots,$ is expressed by formula
 (see \cite[pp. 140, 141]{Kil})
\begin{equation*}
u_n(t)=\varphi_n E_{\alpha, 1} \left(-\lambda_n t^{\alpha}\right)+\psi_n t E_{\alpha, 2} \left(-\lambda_n t^{\alpha}\right)\quad\quad\quad\quad\quad\quad\quad\quad\quad
\end{equation*}
\begin{equation}\label{eq(1.3)}
\quad\quad\quad\quad\quad\quad\quad\quad\quad\quad\quad
+h_n\int\limits_0^t (t-s)^{\alpha-1} E_{\alpha, \alpha} \left(-\lambda_n (t-s)^{\alpha}\right)f(s) d s.
\end{equation}

Thus, the solution to the problem (\ref{eq(10.1)})-(\ref{eq(10.3)}) must be expressed in the form of a formal series
\begin{equation*}
u(t, x)=\sum\limits_{n=1}^{\infty} \Bigg[\varphi_n E_{\alpha, 1} \left(-\lambda_n t^{\alpha}\right)+\psi_n t E_{\alpha, 2} \left(-\lambda_n t^{\alpha}\right)\quad\quad\quad\quad\quad\quad\quad\quad\quad
\end{equation*}
\begin{equation}\label{eq(1.44)}
\quad\quad\quad\quad\quad\quad\quad
+h_n\int\limits_0^t (t-s)^{\alpha-1} E_{\alpha, \alpha} \left(-\lambda_n (t-s)^{\alpha}\right)f(s) d s\Bigg]X_n(x).
\end{equation}

If the series in (\ref{eq(1.44)}) converges uniformly in $\overline{\Omega}_T$, and the series obtained from (\ref{eq(1.44)}) by term-wise application of the operators $\partial_t^{\alpha}$ and $L$ with respect to $x$ also converge uniformly in $\Omega_T$, then the sum of this series belongs to the class $C^{\alpha, 2}_{t, x}\left(\Omega_T\right) \cap C^{1, 0}_{t, x}\left(\overline{\Omega}_T\right)$. Furthermore, it satisfies equation (\ref{eq(10.1)}) as well as the conditions (\ref{eq(10.2)}) and (\ref{eq(10.3)}).

Introduce the notation $D^{k}\left( \cdot\right):=\frac{\partial^{k}\left( \cdot\right)}{\partial x_1^{k_1} \partial x_2^{k_2}\cdots \partial x_N^{k_N}},$ where $ \, \, \sum\limits_{i=1}^Nk_i=k,$ and  $ k_1, \cdots, k_N, k\in \Big\{0, 1, \cdots \Big\}.$

The following statement is true.

\begin{theorem}\label{thm1}
 \it Let  $f(t)\in C[0, T]$ and  the following conditions are met:\\
$
1)\,\,a_{ij}(x)\in C^{\left[\frac{N}{2}\right]+2}\left(\overline{\Omega}\right), \, \, i, j=1,\cdots, N, \, \, c(x)\in C^{\left[\frac{N}{2}\right]+1}\left(\overline{\Omega}\right);
$\\
$
2)\,\,\Big\{ \varphi(x),  \,h(x)\Big\}\in C^{\left[\frac{N}{2}\right]+2}(\overline{\Omega}),\,\, D^{\left[\frac{N}{2}\right]+3}\Big\{ \varphi(x),  \,h(x)\Big\}\in L^2\left(\Omega\right),$ \\
and $\,\Big\{ \varphi(x), \,  \,h(x)\Big\}\Big|_{\partial \Omega}=L\Big\{ \varphi(x), \, \,h(x)\Big\}\Big|_{\partial \Omega} =\cdots=
L^{\left[\frac{N}{4}\right]+1}\Big\{ \varphi(x), \,  \,h(x)\Big\}\Big|_{\partial \Omega}=0.
$\\
$
3)\,\,\psi(x)\in C^{\left[\frac{N}{2}\right]+1}(\overline{\Omega}),\,\, D^{\left[\frac{N}{2}\right]+2}\psi(x)\in L^2\left(\Omega\right),$ \\
and $\,\psi(x)\Big|_{\partial \Omega}=L\psi(x)\Big|_{\partial \Omega} =\cdots=
L^{\left[\frac{N+2}{4}\right]}\psi(x)\Big|_{\partial \Omega}=0.
$

Then, there is a unique classical solution to DP.
\end{theorem}

\textbf{Proof.} Each term of the series (\ref{eq(1.44)}) satisfies equation (\ref{eq(10.1)}) by its construction. Therefore, it suffices to demonstrate the uniform and absolute convergence of the series for $u(t, x)$ and $u_t(t, x)$ in the domain $\overline{\Omega}_T$. Additionally, it is necessary to show the uniform and absolute convergence of the series for $\partial_t^{\alpha}u(t, x)$ and $Lu(t, x)$ in the domain $\Omega_T^{t_0} := \big\{(t, x): t_0 \leq t \leq T, \, x \in \overline{\Omega} \big\}$, where $t_0$ is a sufficiently small positive number. To this end, by formal differentiation of (\ref{eq(1.2)}), we have
\begin{equation}\label{eq(1.22)}
u_t(t, x)=\sum\limits_{n=1}^{\infty} u'_n(t) X_n(x).
\end{equation}
Next,  by acting (formally) on the operators $\partial_t^{\alpha}$ and $L$ to  $(\ref{eq(1.2)})$, in view of
\begin{equation}\label{eq(1.3300)}
\partial_t^{\alpha}u_n(t)=-\lambda_n u_n(t)+h_nf(t)
\end{equation}
which is valid for a solution of (\ref{eq(ef3)}),
 we obtain
\begin{equation}\label{23}
\partial_t^{\alpha}u(x,t)=\sum\limits_{n=1}^{\infty}\partial_t^{\alpha}u(t)X_{n}(x)=-\sum\limits_{n=1}^{\infty}\lambda_n u_n(t)X_{n}(x)+f(t)\sum\limits_{n=1}^{\infty}h_nX_{n}(x),
\end{equation}

\begin{equation}\label{24}
Lu(x,t)=\sum\limits_{n=1}^{\infty}u_{n}(t)LX_{n}(x)=-\sum\limits_{n=1}^{\infty}\lambda_n u_{n}(t)X_{n}(x).
\end{equation}

{\bf Step 1.} First, we prove the convergence of the series $(\ref{eq(1.2)})$, $(\ref{eq(1.22)})$ in $\overline{\Omega}_T$.

In view of (\ref{eq(1.3)}) and Proposition \ref{prop1}, from $(\ref{eq(1.2)})$ it follows,
 \begin{equation}\label{eq(01)}
\left|u(x,t)\right|\leq\sum\limits_{n=1}^{\infty}|u_{n}(t)||X_{n}(x)|\leq\mathcal{K}\sum\limits_{n=1}^{\infty}\big(\left|\varphi_n\right|+\left|\psi_n\right|+\left|h_n\right|\big)|X_{n}(x)|.
 \end{equation}
Here and below the letter $\mathcal{K}$ will denote various positive constants independent on functions $\varphi, \psi, h.$

 By differentiation formulas for expressions containing Mittag-Leffler functions (see, \cite[p. 50]{Kil}), from (\ref{eq(1.44)}),   we find for  (\ref{eq(1.22)})
\begin{equation*}
u_t(t, x)=\sum\limits_{n=1}^{\infty} \Bigg[-\lambda_n\varphi_n t^{\alpha-1} E_{\alpha, \alpha} \left(-\lambda_n t^{\alpha}\right)+\psi_n  E_{\alpha, 1} \left(-\lambda_n t^{\alpha}\right)\quad\quad\quad\quad\quad\quad\quad\quad\quad
\end{equation*}
\begin{equation}\label{eq(1.444)}
\quad\quad\quad\quad\quad
+h_n\int\limits_0^t (t-s)^{\alpha-2} E_{\alpha, \alpha-1} \left(-\lambda_n (t-s)^{\alpha}\right)f(s) d s\Bigg]X_n(x).
\end{equation}
 Based on   Proposition \ref{prop1}, having carried out the estimates in (\ref{eq(1.444)}) we obtain
  \begin{equation}\label{eq(0)}
  \left|u_t(t, x)\right|\leq\mathcal{K}\sum\limits_{n=1}^{\infty}\big(\lambda_n\left|\varphi_n\right|+\left|\psi_n\right|+\left|h_n\right|\big)|X_{n}(x)|.
  \end{equation}
  Under conditions 1)-3) of Theorem \ref{thm1}, the convergence of the series on the right-hand sides of inequalities $(\ref{eq(01)}),$ $(\ref{eq(0)})$ follows from Theorem 8 in \cite{Il30}.  Consequently,  series $(\ref{eq(1.2)})$, $(\ref{eq(1.22)})$  converge uniformly and absolutely in $\overline{\Omega}_T$.  Thus,  $u(t, x)\in C^{1, 0}_{t, x}(\overline{\Omega}_T).$

  \textbf{Step 2.} We will establish the absolute and uniform convergence of the series in (\ref{23}) and (\ref{24}) within the domain $\Omega_T^{t_0}$.
In view of (\ref{eq(1.3)}), (\ref{eq(1.3300)}) and Proposition \ref{prop1}, for the second series in (\ref{eq(1.22)}) we have
 $$
 \left|\partial^{\alpha}_tu(t, x)\right|\leq \mathcal{K}\Bigg[\sum_{n=1}^{\infty}\frac{\lambda_n}{1+\lambda_n t_0^{\alpha}}\left|\varphi_n X_n(x)\right|+\sum_{n=1}^{\infty}\frac{\lambda_n}{1+\lambda_n t_0^{\alpha}}\left|\psi_n X_n(x)\right|
  $$
  \begin{equation}\label{eq(2.3339)}
  +\sum_{n=1}^{\infty}\left|h_n X_n(x)\right|+\sum_{n=1}^{\infty}\lambda_n\left|h_n X_n(x)\right|\int\limits_0^t\frac{s^{\alpha-1}ds}{1+\lambda_ns^{\alpha}}\Bigg]. \end{equation}


Since $\frac{\lambda_n}{1+\lambda_n t^{\alpha}} < t_0^{-\alpha}$, the first three series on the right-hand side of (\ref{eq(2.3339)}) are uniformly convergent due to the proven convergence of the series (\ref{eq(0)}). Applying for the last series in (\ref{eq(2.3339)}) \cite{DurRax}:
\begin{equation*}\label{eq(22.3339)}
\max\limits_{y\ge0}\frac{y^{\theta}}{1+y}=\frac{\left(\frac{\theta}{1-\theta}\right)^{\theta}}{1+\frac{\theta}{1-\theta}},\quad 0<\theta<1,
\end{equation*}
for $\theta=(\alpha-1)/\alpha\in(0, 1)$, we find that this series  is estimated by $\mathcal{K}\sum_{n=1}^{\infty}\lambda^{1/\alpha}_n\left|h_n X_n(x)\right|.$
 Since $\lambda_n^{1/\alpha} < \lambda_n$ for sufficiently large $n \in \mathbb{N}$, the convergence of this series can be established in a manner similar to the convergence of the first series in (\ref{eq(0)}), with the function $\varphi$ replaced by $h$.

Because the series on the right-hand side of (\ref{24}) coincides with the first series on the right-hand side of (\ref{23}), we consider the convergence of the series in (\ref{24}) to be established. Consequently, $u(t, x) \in C^{\alpha, 2}_{t, x}\left(\overline{\Omega}_T\right).$

Due to $u(t, x)\in C^{1, 0}_{t, x}\left(\overline{\Omega}_T\right)$, then by tending $t$ to $0+$ in (\ref{eq(1.44)}) and (\ref{eq(1.444)}) for any $x\in \overline{\Omega},$ we get
\begin{equation*}
\lim\limits_{t\rightarrow 0}u(t, x)=u(0, x)=\sum\limits_{n=1}^{\infty} \varphi_n X_n(x)=\varphi(x),\end{equation*}
\begin{equation*}
\lim\limits_{t\rightarrow 0}u_t(t, x)=u_t(0, x)=\sum\limits_{n=1}^{\infty} \psi_n X_n(x)=\psi(x),
\end{equation*}
that is the initial  conditions   (\ref{eq(10.2)}) are satisfied.

Theorem \ref{thm1} is proved.

{\bf Remark 1.} In the proof of Theorem \ref{thm1}, we have used the estimate in Proposition \ref{prop1} for $z = -\lambda_n t^{\alpha}$. If we require the fulfillment of equation (\ref{eq(10.1)}) in the class of functions $C^{\alpha,\, 2}_{t, \, x}\left(\overline{\Omega}_T\right)$, we must use a rougher estimate, $ \left|E_{\alpha,\beta}\left(-\lambda_n t^{\alpha}\right)\right| \leq C, $ which follows from Proposition \ref{prop1}. This estimate is employed in proving $u(t, x) \in C^{1,\, 0 }_{t,\, x}\left(\overline{\Omega}_T\right)$ (step 1).

An analysis of the proof of Theorem \ref{thm1} indicates that if we replace condition 3) of this theorem with the condition.
\\
$
3')\,\psi(x)\in C^{\left[\frac{N}{2}\right]+2}(\overline{\Omega}),\,\, D^{\left[\frac{N}{2}\right]+3}\psi(x)\in L^2\left(\Omega\right),$
and $\,\psi(x)\Big|_{\partial \Omega}=L\psi(x)\Big|_{\partial \Omega} =\cdots $ \\ $=
L^{\left[\frac{N}{4}\right]+1}\psi(x)\Big|_{\partial \Omega}=0,
$ \\
 Thus, we can conclude that the solution to the DP belongs to the class $u(t, x) \in C^{\alpha,\, 2}_{t, \, x}\left(\overline{\Omega}_T\right).$

Consequently, the following statement derives from Remark 1.
\begin{theorem}\label{thm2}
Let  $f(t)\in C[0, T]$ and  the following conditions are met:\\
$
1)\,\,a_{ij}(x)\in C^{\left[\frac{N}{2}\right]+2}\left(\overline{\Omega}\right), \, \, i, j=1,\cdots, N, \, \, c(x)\in C^{\left[\frac{N}{2}\right]+1}\left(\overline{\Omega}\right);
$\\
$
2)\,\,\Big\{ \varphi(x), \, \psi(x), \,h(x)\Big\}\in C^{\left[\frac{N}{2}\right]+2}(\overline{\Omega}),\,\, D^{\left[\frac{N}{2}\right]+3}\Big\{ \varphi(x), \, \psi(x),  \,h(x)\Big\}\in L^2\left(\Omega\right),$ \\
and $\,\Big\{ \varphi(x), \, \psi(x), \,  \,h(x)\Big\}\Big|_{\partial \Omega}=L\Big\{ \varphi(x), \, \psi(x), \, \,h(x)\Big\}\Big|_{\partial \Omega} =\cdots $ \\ $=
L^{\left[\frac{N}{4}\right]+1}\Big\{ \varphi(x), \, \psi(x), \,  \,h(x)\Big\}\Big|_{\partial \Omega}=0.
$

Then, there is a unique classical solution to DP such that $u(t, x)\in C^{\alpha,\, 2}_{t, \, x}\left(\overline{\Omega}_T\right).$
\end{theorem}

\begin{center}
\textbf{Investigation of IP1}
\end{center}

The following statement is true.
\begin{theorem}\label{thm3}
Let the functions $\varphi(x),$ $\psi(x)$ and $h(x)$ satisfy the conditions of Theorem \ref{thm1}, $h(x) \not \equiv$ 0 on $\overline{\Omega}, \, g(t) \in AC^2[0, T]$ and the conditions
\begin{equation*}\label{eq(2.5)}
\int\limits_{\Omega} \varphi(x) h(x) d x=g(0)=\sum_{n=1}^{\infty}h_n \varphi_n
\end{equation*}
are valid. Then, there exists a unique solution to IP such that $f(t)\in AC[0, T].$
\end{theorem}

{\bf Proof.}
The solution of the direct problem (\ref{eq(10.1)})-(\ref{eq(10.3)}) satisfies the additional condition (\ref{eq(ad1)}). Then, we get \begin{equation}\label{eq(2.1)}
\int\limits_{\Omega} h(x)u(t, x)dx=\int\limits_{\Omega} h(x) \sum\limits_{n=1}^{\infty} u_n(t) X_n(x) d x=\sum_{n=1}^{\infty}h_n u_n(t)=g(t),
\end{equation}
where $h_n=\int_{\Omega} h(x) X_n(x) d x.$

In view of (\ref{eq(1.3)}), and from (\ref{eq(2.1)}), we derive the Volterra integral equation of the first kind with a difference kernel for the unknown function $f(t)$:
\begin{equation}\label{eq(2.2)}
\int\limits_0^t f(s) K(t-s) d s=G(t), \quad t\in [0, T],
\end{equation}
where
\begin{equation}\label{eq(2.3)}
K(t)=t^{\alpha-1}\sum\limits_{n=1}^{\infty} h_n^2  E_{\alpha, \, \alpha}\left(-\lambda_n t^{\alpha}\right),
\end{equation}
 \begin{equation}\label{eq(2.4)}G(t)=g(t)
 -\sum_{n=1}^{\infty}h_n\left[ \varphi_n E_{\alpha, \, 1}\left(-\lambda_n t^{\alpha}\right)+\psi_n \ t E_{\alpha, \, 2}\left(-\lambda_n t^{\alpha}\right)\right].
\end{equation}

Under the conditions imposed on the functions $\varphi(x)$, $\psi(x)$, and $h(x)$ in Theorem \ref{thm1}, the series in (\ref{eq(2.3)}) and (\ref{eq(2.4)}) converge uniformly and allow term-by-term differentiation with respect to $t$. Differentiating (\ref{eq(2.2)}) and using formula (1.10.7) from \cite{Kil}, we obtain
\begin{equation}\label{eq(2.22)}
\int\limits_0^t f(s) K'(t-s) d s=G'(t), \quad t\in [0, T],
\end{equation}
where
\begin{equation*}\label{eq(2.33)}
K'(t)=t^{\alpha-2}\sum\limits_{n=1}^{\infty} h_n^2  E_{\alpha, \, \alpha-1}\left(-\lambda_n t^{\alpha}\right),
\end{equation*}
\begin{equation*}\label{eq(2.455)}G'(t)=g'(t)+\sum_{n=1}^{\infty}h_n\left[ \lambda_n \varphi_nt^{\alpha-1}E_{\alpha, \, \alpha}\left(-\lambda_n t^{\alpha}\right)-\psi_n  E_{\alpha, \, 1}\left(-\lambda_n t^{\alpha}\right)\right].
\end{equation*}

Applying the operator $\partial_t^{\alpha-1}$ to both sides of (\ref{eq(2.22)}) and in view of  Proposition \ref{prop4}, we have
\begin{equation}\label{eq(2.222)}
f(t)\sum\limits_{n=1}^{\infty} h_n^2 +\int\limits_0^t f(s) K_0(t-s) d s=G_0(t),
\end{equation}
where
\begin{equation}\label{eq(2.333)}
K_0(t)=-t^{\alpha-1}\sum\limits_{n=1}^{\infty}\lambda_n h_n^2  E_{\alpha, \, \alpha}\left(-\lambda_n t^{\alpha}\right), \quad t\in [0, T],\end{equation}
$$ G_0(t)=\partial_t^{\alpha-1}G'(t)=\partial_t^{\alpha-1}g'(t)\quad\quad\quad\quad\quad\quad\quad\quad\quad\quad\quad\quad\quad\quad\quad\quad\quad\quad\quad\quad\quad\quad\quad\quad\quad\quad\quad\quad\quad\quad\quad\quad\quad\quad\quad\quad
$$
\begin{equation}\label{eq(2.444)}+\sum_{n=1}^{\infty}h_n\left[- \lambda_n \varphi_nE_{\alpha, \, 1}\left(-\lambda_n t^{\alpha}\right)+\lambda_n\psi_ntE_{\alpha, \, 2}\left(-\lambda_n t^{\alpha}\right)\right],  t\in [0, T].
\end{equation}

In the previous relations, the formulas for fractional differentiation from Proposition \ref{prop3} were applied. By Theorem \ref{thm2}, the function $h(x) \neq 0$ on $\overline{\Omega}$. Therefore, by the Parseval-Steklov equality, we get
\begin{equation*}
\sum_{n=1}^{\infty} h_n^2=\int_{\Omega} h^2(x) d x=\|h\|_{L_2(\Omega)}^2>0.
\end{equation*}

Let us demonstrate that the series in (\ref{eq(2.333)}) and (\ref{eq(2.444)}) converge absolutely and uniformly, thereby representing continuous functions. For the series in (\ref{eq(2.333)}), using Proposition \ref{prop1}, we obtain
\begin{equation*}
\left|\sum\limits_{n=1}^{\infty}\lambda_n h_n^2 E_{\alpha, \, \alpha}\left(-\lambda_n t^{\alpha}\right)\right|\leq C \sum\limits_{n=1}^{\infty}\lambda_n h_n^2, \quad t\in[0, T].
\end{equation*}

To prove the convergence of the  numerical series $\sum\limits_{n=1}^{\infty}\lambda_n h_n^2$, consider the following functional:
\begin{equation*}
J(Y)=\int\limits_{\Omega}\Big[\left|\mathbf{A}^{\frac{1}{2}}(x)\nabla Y(x)\right|^2+c(x)Y^2(x)\Big]dx.
\end{equation*}
Substitute the expression $Y(x) = h(x) - \sum\limits_{j=1}^{n-1} h_j X_j(x)$ into it, and use the identities 
$$
\Big(\text{div}\left[\mathbf{A}(x) \nabla X_n(x)\right], X_n(x)\Big) 
$$$$
= \Big(\mathbf{A}^{1/2}\mathbf{A}^{1/2} \nabla X_n(x), \nabla X_n(x)\Big) = \Big(\mathbf{A}^{1/2} \nabla X_n(x), \mathbf{A}^{1/2}\nabla X_n(x)\Big) = \Big\| \left|\mathbf{A}^{1/2} \nabla X_n\right| \Big\|_{L_2(\Omega)},
$$ 
which follow from (\ref{5}), the symmetry $\mathbf{A} = \mathbf{A}^{\top}
$, 
and the positive definiteness of $\mathbf{A}(x)$. From this and (\ref{5}), it follows that
$$\lambda_n=\Big\| \ \left|\mathbf{A}^{\frac{1}{2}} \nabla X_n\right| \ \Big\|_{L_2(\Omega)}+\Big\| \ \left|c^{\frac{1}{2}} \nabla X_n\right| \ \Big\|_{L_2(\Omega)}.$$
  Based on  the easily verified  identity
\begin{equation*}\Big(\mathbf{A}(x) \nabla X_k(x), \ \nabla X_j(x)\Big)+\Big(c(x)X_k(x), \ X_j(x)\Big)=0, \, \, \, k\neq j, \end{equation*}
we get\begin{equation*}\label{eq(3.01)}
J\left[h(x)-\sum\limits_{j=1}^{n-1}h_jX_j(x)\right]=\Big\| \ \left|\mathbf{A}^{\frac{1}{2}} \nabla h\right| \ \Big\|_{L_2(\Omega)}+\Big\| \ \left|c^{\frac{1}{2}} \nabla h\right| \ \Big\|_{L_2(\Omega)}
\end{equation*}
\begin{equation*}\label{eq(3.01)}
+\sum\limits_{j=1}^{n-1}\lambda_jh_j^2+2\sum\limits_{j=1}^{n-1}h_j\Bigg[\Big(\mathbf{A}(x)\nabla h(x), \  \nabla X_j(x) \Big)+\Big(c(x)h(x), \ X_j(x)\Big)\Bigg].
\end{equation*}
Continuing the calculations and using the conditions of Theorem \ref{thm1} for $h(x)$, as well as the equality (\ref{5}) with $X = X_n(x)$ and $\lambda = \lambda_n$, we obtain
\begin{equation}\label{eq(3.03)}
J\left[h(x)-\sum\limits_{j=1}^{n-1}h_jX_j(x)\right]=\Big\| \ \left|\mathbf{A}^{\frac{1}{2}} \nabla h\right| \ \Big\|_{L_2(\Omega)}+\Big\| \ \left|c^{\frac{1}{2}} \nabla h\right| \ \Big\|_{L_2(\Omega)}
-\sum\limits_{j=1}^{n-1}\lambda_jh_j^2.
\end{equation}
Since $J(Y) \geq 0$, it follows from formula (\ref{eq(3.03)}) the inequality
\begin{equation}\label{eq(11.011)}
\sum\limits_{j=1}^{\infty}\lambda_jh_j^2\leq\Big\| \ \left|\mathbf{A}^{\frac{1}{2}} \nabla h\right| \ \Big\|_{L_2(\Omega)}+\Big\| \ \left|c^{\frac{1}{2}} \nabla h\right| \ \Big\|_{L_2(\Omega)},
\end{equation}
and the convergence of the series on the left side. The relation (\ref{eq(11.011)}) is a Bessel-type inequality for the function $h(x)$ (see also Lemma 5 in \cite{Il30}).

To prove the uniform convergence of the series in (\ref{eq(2.444)}), it suffices, as follows from Proposition \ref{prop1}, to establish the convergence of the same type of majorant numerical series $C \sum_{n=1}^{\infty} \lambda_n |h_n| |\varphi_n|$ and $C \sum_{n=1}^{\infty} \lambda_n |h_n| |\psi_n|$. The convergence of these numerical series can be demonstrated as follows.
Using, for instance, the Cauchy inequality for the first series, we have
\begin{equation*}\label{eq(2.44)}
C\sum\limits_{n=1}^{\infty}\lambda_n|h_n|  |\varphi_n|=C\sum\limits_{n=1}^{\infty}\sqrt{\lambda_n}|h_n|  \sqrt{\lambda_n}|\varphi_n|\leq C\sqrt{\sum\limits_{n=1}^{\infty}\lambda_nh_n^2}\, \,   \sqrt{\sum\limits_{n=1}^{\infty}\lambda_n\varphi_n^2}.
\end{equation*}
The convergence of the series on the right-hand side of this relation is directly ensured by the preceding reasoning.

Since $g(t) \in AC^2[0, T]$, it follows from Lemma 2 that $\partial_t^{\alpha-1} g(t) \in AC[0, T]$. The absolute continuity of the sum of the uniformly convergent series in (\ref{eq(2.333)}) and (\ref{eq(2.444)}) is evident. Hence, equation (\ref{eq(2.222)}) is an integral equation of the second kind of Volterra type, with absolutely continuous input data $K_0(t)$ and $G_0(t)$. Such an equation possesses a unique solution within the class of absolutely continuous functions, and the solution is given by the formula
\begin{equation*}\label{eq(2.6)}
f(t)=\frac{G_0(t)}{\|h\|_{L_2(\Omega)}^2}+\int_0^tR(t-s) \frac{G_0(s)}{\|h\|_{L_2(\Omega)}^2} d s,
\end{equation*}
where $R(t)$ is the resolvent of the function $-\frac{K_0(t)}{\|h\|_{L_2(\Omega)}^2}.$  Theorem \ref{thm2} is proved.

\begin{center}
 \textbf{Uniqueness of IP1 solution}
\end{center}
\begin{theorem}\label{thm4}
Let  $h(x)\neq 0$ on $\overline{\Omega}.$
Then, if a solution to IP1 exists in the class $\left\{u(t, x), \, f(t)\right\}\in C^{\alpha,\, 2}_{t, \, x}\left(\overline{\Omega}_T\right)\times C[0, T]$, it is unique within the this class of functions.
\end{theorem}

To prove the uniqueness of the solution to the IP1, we assume, for the sake of contradiction, that there exist two distinct solutions  $\left\{u_1(t, x), \, f_1(t)\right\}$ and $\left\{ u_2(t, x), \, f_2(t)\right\}$, where  $u_j(t, x), \, j=1, 2,$  are the corresponding to $f_i,$  solutions to DP (\ref{eq(10.1)})-(\ref{eq(10.3)}), and that these solutions also satisfy the integral conditions $$\int_{\Omega} h(x)(u_j)_t(t, x)dx=g(t), \, \, j=1, 2, \, \,   t\in [0, T].$$
Then,  the function
$u(t, x) = u_1(t, x)-u_2(t, x)$ solves the equation
(\ref{eq(10.1)}), where    $f(t)=f_1(t)-f_2(t)$
 subject to the homogeneous boundary conditions  (\ref{eq(10.2)}) $\left(\varphi(x)\equiv0, \ \psi(x)\equiv0\right)$, (\ref{eq(10.3)}) and also satisfies  the integral condition $$\int_{\Omega} h(x)u(t, x)dx=0.$$

We define the functions  $u_n(t)=\int_{\Omega} u(t, x) X_n(x) d x,$ where $X_n(x)$ are eigenfunctions of the spectral problem $(\ref{5}).$ By applying the fractional operator $\partial^{\alpha}_t$ to $u_n(t),$ and using equation (\ref{eq(10.1)}), we obtain the following result:
\begin{equation}\label{uniq1}
\partial^{\alpha}_t u_n(t)=\int_{\Omega}\Big[ L u(t, x)+f(t)h(x)\Big] X_n(x) d x=-\lambda_n u_n(t)+h_nf(t).
\end{equation}
Let us consider these equations with homogeneous initial conditions
\begin{equation}\label{uniq2}
 u_n(0)=u'_n(0)=0,
\end{equation}
and overdetermination condition
\begin{equation}\label{uniq3}
\int_{\Omega} h(x)u(t, x)dx=0.
\end{equation}
The solution of the problem (\ref{uniq1}), (\ref{uniq2}) is expressed by the formula
\begin{equation}\label{uniq4}
u_n(t)=h_n\int\limits_0^t (t-s)^{\alpha-1} E_{\alpha, \alpha} \left(-\lambda_n (t-s)^{\alpha}\right)f(s) d s.
\end{equation}
By using (\ref{uniq3}) and the relationship derived from it in (\ref{eq(2.1)}) with $g(t)\equiv 0,$ we can simplify the last equation into the form of the homogeneous equation (\ref{eq(2.2)}), which represents the Volterra integral equation of the first kind with a difference kernel with respect to $f(t).$ After applying the transformations outlined in Section 3, this equation is further reduced to a homogeneous Volterra integral equation of the second kind, corresponding to (\ref{eq(2.222)}). As is well-known, such an equation admits only the trivial solution $f(t)\equiv 0.$ Therefore, from (\ref{uniq4}) we conclude that  $u_n(t)\equiv 0$. Since $X_{n}(x),$ for $n\in \mathbb{N},$ forms a complete orthonormal system in  $L_2(\Omega),$ it follows that  $u(x,t)\equiv0$ in $\Omega_T.$ Thus Theorem \ref{thm3} is proven.

\begin{center}
\textbf{Uniqueness of IP2 solution}
\end{center}
\begin{theorem}\label{thm5}
  Let  $f(t)\neq 0$ on $[0, T].$
Then, if a solution to IP2 exists in the class $\left\{u(t, x), \, h(x)\right\}\in C^{\alpha,\, 2}_{t, \, x}\left(\overline{\Omega}_T\right)\cup C\left(\overline{\Omega}_T\right)$, it is unique.
\end{theorem}

\textbf{Proof.} Here we also assume, to the contrary, that there exist two
solutions $\left\{u_1(t, x), \, h_1(x)\right\}$ and $\left\{ u_2(t, x), \, h_2(x)\right\}$, where  $u_j(t, x), \, j=1, 2,$  are the corresponding to $h_i,$  solutions to DP (\ref{eq(10.1)})-(\ref{eq(10.3)}), and they also satisfy the integral conditions:
$$\int_0^T f(t)(u_j)_t(t, x)dt=\omega(x), \, \, j=1, 2, \, \,   x\in \Omega.$$
Then,  the function
$u(t, x) = u_1(t, x)-u_2(t, x)$ is  the solution to equatin
(\ref{eq(10.1)}), where    $h(x)=h_1(x)-h_2(x)$
  subject to the homogeneous boundary conditions  (\ref{eq(10.2)}) $\left(\varphi(x)\equiv0, \ \psi(x)\equiv0\right)$ and (\ref{eq(10.3)}).  Furthermore, $u$ satisfies the integral condition $$\int\limits_0^T f(t)u_t(t, x)dt=0.$$

Multiply both sides of equation (\ref{eq(10.1)})  by $2u_t(t,x)$ where $u$ is a solution to DP. Then, integrate the resulting expression with respect to  $(t, x)$ over $\Omega_T.$ This leads to the integral identity:
\begin{equation}\label{eq(10.5)}
 2\int\limits_{\Omega}\int\limits_0^Tu_t\, \partial_t^{\alpha}u dtdx+2\int\limits_{\Omega}\int\limits_0^T\, (-Lu)\,u_tdtdx
  =2\int\limits_{\Omega}h(x)\left(\int\limits_0^Tf(t)u_t(t, x)dt\right)dx.
 \end{equation}
Since $u(t, x)\in C^{\alpha,\, 2}_{t, \, x}\left(\overline{\Omega}_T\right)$ (see Theorem 2), all integrals in (\ref{eq(10.5)}) are well defined.
Here,  involving  the following identities
  \begin{equation}\label{eq(10.6)}
  2\int\limits_{\Omega}\int\limits_0^T\, \Big[div \left(-\mathbf{A}(x) \ \nabla u(t, x)\right)u_t+c(x)u \, u_t \Big]dtdx
 \end{equation}
 \begin{equation*}
=2\int\limits_{\Omega}\int\limits_0^T\Big[\mathbf{A}(x) \nabla u\Big] \nabla u_tdtdx +\int\limits_{\Omega}\int\limits_0^T\Big(c(x)u^2\Big)_tdtdx,
  \end{equation*}

\begin{equation*}\label{eq(10.7)}
  \int\limits_0^Tf(t)u_t(t, x)dt=0,
   \end{equation*}
in view of $\left[\mathbf{A}(x) \nabla u\right] \nabla u_t=\left(\mathbf{A}^{1/2}\mathbf{A}^{1/2} \nabla u\right) \nabla u_t=\left(\mathbf{A}^{1/2} \nabla u\right) \left(\mathbf{A}^{1/2}\nabla u_t\right)=
\Big[\left(\mathbf{A}^{1/2} \nabla u\right)^2\Big]_t,$ we can obtain from (\ref{eq(10.5)})
  \begin{equation}\label{eq(10.66)}
 2\int\limits_{\Omega}\int\limits_0^Tu_t\, \partial_t^{\alpha} u dtdx+ \int\limits_{\Omega}\int\limits_0^T\Big[\left(\mathbf{A}^{\frac{1}{2}}(x)\nabla u\right)^2 \Big]_t dtdx
  +\int\limits_{\Omega}\int\limits_0^T\left(c(x)u^2\right)_tdtdx=0.
  \end{equation}

  Note that the matrix $\mathbf{A}(x)$
 is positive definite in $\Omega$, thus there exists a unique positive definite square root
$\mathbf{A}^{1/2}$.

Attracting (\ref{eq(1.4444)}) for $u(t,x)$ and integrating it over $\Omega\times(0, T)$,
 we get
 \begin{equation}\label{eq(10.8)}
 \int\limits_{\Omega}\int\limits_0^T u_t(t,x)\partial_t^{\alpha}u(t,x)dtdx\geq \frac{1}{2}\int\limits_{\Omega}\int\limits_0^T\partial_t^{\gamma}u_t^2(t,x)dtdx, \ \ \gamma=\alpha-1. \end{equation}
Further from  (\ref{eq(10.66)}) in view of (\ref{eq(10.88)}),  the homogeneous initial, boundary  conditions,  we arrive at  the inequality
\begin{equation}\label{eq(10.88)}  \int\limits_{\Omega}\int\limits_0^T\partial_t^{\alpha-1}u_t^2(t,x)dtdx+\int\limits_{\Omega}\Big[\left(\mathbf{A}^{\frac{1}{2}}(x) \nabla u(T, x)\right)^2 +c(x)u^2(T, x)\Big]dx\leq 0.
\end{equation}

By the estimate  (\ref{eq(41.444)}) (Proposition \ref{prop5}) taking $u_t(0, x)=0$ into consideration, we finally have
\begin{equation*}\label{eq(11.88)}
\frac{1}{2\Gamma(2-\alpha)}\int\limits_{\Omega}\int\limits_0^T\frac{u_t^2(t, x)}{(T-t)^{\alpha-1}}dt dx\end{equation*}
\begin{equation}\label{eq(11.88)}
+\int\limits_{\Omega}\Big[\left(\mathbf{A}^{\frac{1}{2}}(x) \nabla u(T, x)\right)^2 +c(x)u^2(T, x)\Big]dx\leq 0.
\end{equation}

This inequality implies   that $u(t, x)\equiv 0$ in $\overline{\Omega}_T$, and from equation (\ref{eq(10.1)})   under the condition $f(t)\neq 0$ on
$[0, T],$  we have $h(x)\equiv 0.$

{\bf Remark 2.} In previous works \cite{Slod}, \cite{Rom}, the method of proving the uniqueness of the IP2 solution was employed to study the uniqueness of the inverse problem of identifying a space-dependent source for the second-order hyperbolic equation with a damping term, based on the final-time overdetermination condition. It was demonstrated that this method is not applicable in the absence of the damping term. The overdetermination condition used in IP2 for $f=const$ leads to $u(T, x)=\omega(x)+\varphi(x)$. Upon further examination of the uniqueness for IP2, it becomes evident that the damping term is not necessary in the fractional wave equation (1).

For a solution of the hyperbolic equation corresponding to  (1), we have the following limit:
$$\lim\limits_{\alpha\rightarrow 2-0}\int\limits_0^T\frac{1}{\Gamma(2-\alpha)}\frac{u_t^2(t, x)}{(T-t)^{\alpha-1}}dt=u_t^2(T, x),$$ since $$\lim\limits_{\alpha\rightarrow 2-0}\frac{1}{\Gamma(2-\alpha)}\frac{1}{(T-t)^{\alpha-1}}=\delta(T-t)$$ in the generalized sense, $\delta(\cdot)$ denotes the Dirac delta function  \cite{Gel}.
Given that  $u(T, x)=0$, we can deduce from  (\ref{eq(11.88)}) that  $\Big\|u_t(T)\Big\|_{L^2(\Omega)}=0.$
However, this does not imply that $u=0,$  as it only ensures that the time derivative of $u$ vanishes at
$T,$ not the function itself.

\begin{center}
\textbf{Conclusion}
\end{center}

In this paper, we have studied the inverse problems of determining the time-dependent source function (\( t \)-source problem) and the space-dependent source function (\( x \)-source problem) for a time-fractional diffusion-wave equation. These problems were considered under special integral overdetermination conditions, specifically spatial-averaged and time-averaged integral constraints. We established theorems on the existence and uniqueness of solutions for IP1, and the uniqueness result for IP2.

A natural extension of this work would be to investigate the existence and uniqueness of solutions for these inverse problems when considering a more general operator \( L \) with coefficients that also depend on the time variable. However, the proof techniques developed in this paper are not directly applicable in such a setting.

For subdiffusion equations where the coefficients depend on both spatial and temporal variables, interesting uniqueness results have been established for determining an unknown space-dependent source function under a final-time overdetermination condition, as demonstrated in \cite{Kil12, Kil122}. Exploring the existence of solution of  similar inverse problems    for the time-fractional diffusion-wave equation with time-dependent coefficients remains an open and challenging problem for future research.

\end{document}